\documentclass[12pt,reqno]{amsart}
 \usepackage{graphics}
\usepackage{amsfonts, amsmath,enumerate}
\usepackage{amsrefs}
\usepackage{mathtools}

\numberwithin{equation}{section}
\newtheorem{theorem}{Theorem}[section]
\newtheorem{proposition}[theorem]{Proposition}

\newtheorem*{remarks}{Remarks}
\newtheorem*{notes}{Notes \& Comments}

\def\a{\mathfrak{a}}

\def\grad{{\rm{grad}}}
\def\Ric{{\rm{Ric}}}

\def\n{\nabla}
\def\scal{\rm{scal}}
\def\<{\langle}
\def\>{\rangle}

\def\a{\alpha}

\def\bea{\begin{eqnarray*} }
\def\eea{\end{eqnarray*} }

\def\beq{\begin{equation}}

\def\H{\mathord{\mathbb H}}

\def\Q{\mathord{\mathbb Q}}
\def\R{\mathord{\mathbb R}}

\def\S{\mathord{\mathbb S}}

\def\2{\frac{1}{2}}
\def\3{{\ss}}

\def\.{\cdot}

\def\<{\langle}
\def\>{\rangle}

\def\be{\begin{equation}}
\def\ee{\end{equation}}
\def\bea{\begin{eqnarray}}
\def\eea{\end{eqnarray}}
\def\bsm{\left(\begin{smallmatrix}}
\def\esm{\end{smallmatrix}\right)}
\def\bsmp{\left[\begin{smallmatrix}}
\def\esmp{\end{smallmatrix}\right]}
\def\bpm{\begin{pmatrix}}
\def\epm{\end{pmatrix}}

\newcommand{\quotes}[1]{``#1''}

\begin{document}

\title[On Complete conformally flat submanifolds with nullity]{On Complete Conformally flat submanifolds with nullity in Euclidean space}
\author{Christos-Raent Onti}
\date{}
\maketitle

\begin{abstract}
In this note, we investigate conformally flat submanifolds of Euclidean space with positive 
index of relative nullity. Let $M^n$ be a complete conformally flat manifold and let 
$f\colon M^n\to \R^m$ be an isometric immersion. We prove the following results: 
(1) If the index of relative nullity is at least two, then $M^n$ is flat and $f$ is a cylinder over a 
flat submanifold. 
(2) If the scalar curvature of $M^n$ is non-negative and the index of relative nullity
  is positive, then $f$ is  a cylinder over a submanifold with constant non-negative sectional curvature.
(3) If the scalar curvature of $M^n$ is non-zero
and the index of relative nullity is constant and equal to one, then $f$ is a cylinder over a 
$(n-1)$-dimensional submanifold with non-zero constant sectional curvature. 
\end{abstract}

\renewcommand{\thefootnote}{\fnsymbol{footnote}} 
\renewcommand{\thefootnote}{\arabic{footnote}} 

\renewcommand{\thefootnote}{\fnsymbol{footnote}} 
\renewcommand{\thefootnote}{\arabic{footnote}}

\section{Introduction}

A Riemannian manifold $M^n$ is said to be \emph{conformally flat} if each 
point lies in an open neighborhood conformal to an open subset of the Euclidean 
space $\R^n$. The geometry and topology of such Riemannian manifolds 
have been investigated by several authors from the intrinsic point of view. Some  
of the many papers are  \cite{cat16,ku49,ku50,cadjnd11,cahe06,no93,sy88}. 

\vspace{1ex}

Around 1919, Cartan \cite{car17} initiated the investigation of such Riemannian manifolds 
from the submanifold point of view by studying the case of conformally flat 
Euclidean hypersurfaces (see also \cite{mmf85,pin85}). In 1977, Moore \cite{mo77} 
extended Cartan's result in higher (but still low) codimension (see also \cite{df96,df99,mm78}). 
Recently, the author, in collaboration with Dajczer and Vlachos, investigated in \cite{dov18} the case of
conformally flat submanifolds with flat normal bundle in arbitrary codimension (see also \cite{dote11}). 

\vspace{1ex}

In this short note, we address and deal with the following:
\vspace{1ex}

{\noindent{\bf Problem.}}
{\it Classify complete conformally flat submanifolds of Euclidean space 
with positive index of relative nullity and arbitrary codimension.
}
\vspace{1ex}

Recall that the {\it index of relative nullity} at a point $x\in M^n$ of a submanifold $f\colon M^n\to \R^m$ 
is defined as the dimension of the kernel of its second fundamental form $\alpha\colon TM\times TM\to N_fM$, with values in 
the normal bundle.
\vspace{1ex}

The first result provides a complete answer in the case where 
the index of relative nullity is at least two, and is stated as follows:

\begin{theorem}\label{main2}
Let $M^n$ be a complete, conformally flat manifold and let $f\colon M^n\to\R^m$ be an 
isometric immersion with index of relative nullity at least two at any 
point of $M^n$.  Then $M^n$ is flat and $f$ is a cylinder over a flat submanifold.
\end{theorem}

The next result provides a complete answer in the case where the scalar curvature is non-negative.

\begin{theorem}\label{main1}
Let $M^n$ be a complete, conformally flat manifold 
with non-negative scalar curvature and let $f\colon M^n\to\R^m$ be an 
isometric immersion with positive index of relative nullity.  Then 
$f$ is a cylinder over a submanifold with constant non-negative sectional curvature. 
\end{theorem}

Observe that there are complete conformally flat manifolds such that the scalar curvature is non-negative 
while the sectional curvature is not. Easy examples are the Riemannian products
$M^n=\S^{n-m}\times\H^m,\ n\geq 2m,$ where $\S^{n-m}$ and $\H^m$ are the sphere and the 
hyperbolic space of sectional curvature $1$ and $-1$, respectively.

\vspace{1ex}

Finally, the next result provides a complete answer (both local and global) in the case 
where the scalar is non-zero and the index of relative nullity is constant and equal to one.

\begin{theorem}\label{main3}
Let $M^n$ be a conformally flat manifold with non-zero scalar curvature and let $f\colon M^n\to\R^m$ be an 
isometric immersion with constant index of relative nullity equal to one. Then $f$ is locally either a cylinder 
over a $(n-1)$-dimensional submanifold with non-zero constant sectional curvature or a cone over a 
$(n-1)$-dimensional spherical submanifold with constant sectional curvature. Moreover, if $M^n$ is complete, then $f$ is globally 
a cylinder over a $(n-1)$-dimensional submanifold with non-zero constant sectional curvature.
\end{theorem}

\begin{remarks}{\rm
(I)  If the ambient space form in Theorem \ref{main1} is 
replaced by the sphere $\S_c^m$ of constant sectional curvature $c$, 
then an intrinsic classification can be obtained, 
provided that ${\scal}(M^n)\geq c(n-1)$. This classification follows from a result of Carron and Herzlich 
\cite{cahe06}, since in this case $M^n$ turns out to have non-negative Ricci curvature.
However, we do not obtain any (direct) information on the immersion $f$. \\[1mm]
(II) If $f\colon M^n\to\Q_c^m$ is an isometric immersion of a conformally flat manifold into a 
space form of constant sectional curvature $c$, then one can prove the following: 
(i) If the index of relative nullity is at least two, then $M^n$ has constant sectional curvature $c$. 
In particular, if $f$ is also minimal, then $f$ is totally geodesic. 
(ii) If the index of relative nullity is constant and equal to one, then $f$ is a $1$-generalized cone 
(for the definition, see \cite{dov18}) over an isometric immersion $F\colon \Sigma^{n-1}\to \Q_{\tilde c}^{m-1}$ 
into an umbilical submanifold of $\Q_c^m$. }
\end{remarks}

\begin{notes}
{\rm
The special case of minimal conformally flat hypersurfaces 
$f\colon M^n\to \Q_c^{n+1},$ $n\geq 4,$ 
was treated by do Carmo and Dajczer in \cite{DCD} (without any additional assumption on the index of relative nullity), 
where they showed that these are actually generalized 
catenoids, extending that way a previous result due to Blair \cite{Blair} for the case $c=0$. 

For the \quotes{neighbor} class of Einstein manifolds one can prove that: 
any minimal isometric immersion $f\colon M^n\to \Q_c^m$ of an 
Einstein manifold with positive index of relative nullity is totally geodesic. 
A related result of Di Scala \cite{discala}, in the case where the ambient space 
is the Euclidean one, states that: any minimal isometric immersion $f\colon M^n\to \R^m$ 
of a K\"{a}hler-Einstein manifold is totally geodesic. However, it is not yet known if the 
assumption on K\"{a}hler can be dropped (this was conjectured by Di Scala in the same paper). 
Of course, in some special cases the conjecture is true, as have already been pointed out in \cite{discala}.
Finally, we note that Di Scala's theorem still holds true if the K\"{a}hler (intrinsic) assumption 
is replaced by the (extrinsic) assumption on $f$ having flat normal bundle. This follows directly from N\"{o}lker's 
theorem \cite{no90}, since, in this case, $f$ has homothetical Gauss map.
}
\end{notes}

\section{Preliminaries} 

In this section we recall some basic facts and definitions. Let $M^n$ be a Riemannian manifold
and let $f\colon M^n\to \R^m$ be an isometric immersion. The \emph{index of relative nullity} 
$\nu(x)$ at $x\in M^n$ is the dimension of the \emph{relative nullity subspace} 
$\Delta(x)\subset T_xM$ given by
$$
\Delta(x)=\{X\in T_xM: \alpha(X,Y)=0\;\;\mbox{for all}\;\;Y\in T_xM\}.
$$
It is a standard fact that on any open 
subset where the index of relative nullity $\nu(x)$ is constant, the relative nullity distribution 
$x\mapsto \Delta(x)$ is integrable and its leaves are totally geodesic in $M^n$ and $\R^m$. 
Moreover, if $M^n$ is complete then the leaves are also complete along the open subset 
where the index reaches its minimum (see \cite{dajczer}).
If $M^n$ splits as a Riemannian product $M^n=\Sigma^{n-k}\times \R^k$ and there is an isometric 
immersion $F\colon \Sigma^{n-k}\to \R^{m-k}$ such that $f=F\times {\rm id}_{\R^k}$, then we say that 
$f$ is a $k$-cylinder (or simply a cylinder) over $F$.
 
\vspace{1ex}

The following is due to Hartman \cite{har70}; cf. \cite{dt}.

\begin{theorem}\label{hartman}
Let $M^n$ be a complete Riemannian manifold with non-negative Ricci curvature and let 
$f\colon M^n\to \R^m$ be an isometric immersion with minimal index of relative nullity $\nu_0>0$.
Then $f$ is a $\nu_0$-cylinder.
\end{theorem}

A smooth tangent distribution $D$ is called {\it totally umbilical} if there exists a smooth section 
$\delta\in \Gamma(D^\perp)$ such that 
$$
\<\n_X Y, T\>=\<X,Y\>\<\delta,T\>
$$  
for all $X,Y\in D$ and $T\in D^\perp$. The following is contained in \cite{dt}.

\begin{proposition}\label{prop}
Let $f\colon M^n\to \R^m$ be an isometric immersion of a Riemannian manifold 
with constant index of relative nullity $\nu=1$. Assume that the conullity distribution 
$\Delta^\perp$ is totally umbilical (respectively, totally geodesic). Then $f$ is locally
a cone over an isometric immersion $F\colon \Sigma^{n-1}\to\S^{m-1}\subset\R^m$ 
(respectively, a cylinder over an isometric immersion $F\colon \Sigma^{n-1}\to\R^{m-1}\subset\R^m$).
\end{proposition}

We also need the following two well-known results; cf. \cite{dt}.

\begin{proposition}\label{conflat}
A Riemannian product is conformally flat if and only if one of the following 
possibilities holds:
\begin{enumerate}[(i)]
  \item One of the factors is one-dimensional and the other one has constant sectional curvature.
  \item Both factors have dimension greater than one and are either both flat or have opposite 
  constant sectional curvatures.
\end{enumerate}
\end{proposition}

\begin{proposition}\label{conflat2}
Let $M=M_1\times_\rho M_2$ be a warped product manifold. If $M_1$ has dimension one 
then $M$ is conformally flat if and only if $M_2$ has constant sectional curvature.
\end{proposition}

\section{The proofs} 

Let $M^n$ be a conformally flat manifold and let $f\colon M^n\to \R^m$ be an isometric immersion. 
It is well-known that in this case the curvature tensor has the form
$$
R(X,Y,Z,W) = L(X,W)\<Y,Z\>-L(X,Z)\<Y,W\>+L(Y,Z)\<X,W\>-L(Y,W)\<X,Z\>  
$$
in terms of the Schouten tensor given by
\be\label{shouten}
L(X,Y)=\frac{1}{n-2}\left(\Ric(X,Y)-\frac{s}{2(n-1)}\<X,Y\>\right)
\ee
where $s$ denotes the scalar curvature. In particular, 
the sectional curvature is given by
\be\label{seccur}
K(X,Y)=L(X,X)+L(Y,Y)
\ee
where  $X,Y\in TM$ are orthonormal vectors. 

A straightforward computation of the Ricci tensor using the 
Gauss equation
\be\label{eqgauss}
R(X,Y,Z,W)=\<\alpha(X,W),\alpha(Y,Z)\>-\<\alpha(X,Z),\alpha(Y,W)\> 
\ee
yields
\be\label{ricci}
\Ric(X,Y)
=\<nH,\alpha(X,Y)\>-\sum_{j=1}^n\<\alpha(X,X_j),\alpha(Y,X_j)\>
\ee
where $X_1,\dots,X_n$ is an 
orthonormal tangent basis. 

We obtain from \eqref{seccur} and \eqref{eqgauss} that
\be\label{dec2}
L(X,X)+L(Y,Y)=\<\a(X,X),\a(Y,Y)\>-\|\a(X,Y)\|^2
\ee
for any pair $X,Y\in TM$ of orthonormal vectors. Using \eqref{shouten} it follows from \eqref{dec2} that
\be
\Ric(X,X)+\Ric(Y,Y) = \frac{s}{n-1}+(n-2)(\<\a(X,X),\a(Y,Y)\>-\|\a(X,Y)\|^2)  \label{ric}
\ee
for any pair $X,Y\in TM$ of orthonormal vectors. Now, assume that $\nu>0$ and choose a 
unit length $X\in \Delta$. Using \eqref{ricci}, it follows from \eqref{ric} that 
\be\label{ric1}
\Ric(Y,Y) = \frac{s}{n-1}
\ee
for all unit length $Y\perp X$.

\vspace{1ex}
{\noindent {\it Proof of Theorem \ref{main2}:}}
It follows from \eqref{ricci} and \eqref{ric1} that $s=0$. Thus, it follows from \eqref{ric1} that 
$M^n$ is Ricci flat. Since $M^n$ is conformally flat we obtain that $M^n$ is flat. The desired result follows from Theorem 
\ref{hartman} and Proposition \ref{conflat}. \qed

\vspace{1ex}

{\noindent {\it Proof of Theorem \ref{main1}:}} 
It follows from \eqref{ric1} that $\Ric\geq 0$. The desired result follows from 
Theorem \ref{hartman} and Proposition \ref{conflat}. \qed

\vspace{1ex}
{\noindent {\it Proof of Theorem \ref{main3}:}}
It follows from \eqref{dec2}, \eqref{shouten} and \eqref{ricci} that 
\be\label{eq1}
L(Y,Y)=-L(X,X)=\frac{s}{2(n-1)(n-2)}=:h\ \ \text{and} \ \ L(X,Y)=0
\ee
for any unit length vectors $X\in \Delta$ and $Y\in \Delta^\perp$. Moreover, we have that
\be\label{eq2}
L(Y,Z)=0
\ee 
for any pair $Y,Z\in \Delta^\perp$ of orthonormal vectors. Indeed, if $Y$ and $Z$ are two such 
vectors then using \eqref{eq1} we get
$$
L(Y+Z,Y+Z)=2h,
$$
and \eqref{eq2} follows.
Now, since $M^n$ is conformally flat we have that 
$L$ is a Codazzi tensor. Thus
$$
(\n_Y L)(X,X)=(\n_X L)(Y,X)
$$
for all unit length $X\in \Delta$ and $Y\in \Delta^\perp$. It follows, using \eqref{eq1} and \eqref{eq2} that
$$
Y(h)=0,
$$
for all unit length $Y\in \Delta^\perp$. Therefore
$$
h=h(t)=h(\gamma(t)),\ t\in I,
$$ 
where $\gamma\colon I\subset \R\to M^n$ 
is a leaf of the nullity distribution $\Delta$, parametrized by arc length. Using again the fact that $L$ is a Codazzi tensor, 
we get 
$$
(\n_{\gamma'(t)} L)(Y,Z)=(\n_Y L)(\gamma'(t),Z)
$$
for all unit length $Y,Z\in \Delta^\perp$, or equivalently, 
$$
\<Y,Z\>\<\grad\log\sqrt{|h(t)|},\gamma'(t)\>=\<\n_Y Z,\gamma'(t)\>
$$
for all unit length $Y,Z\in \Delta^\perp$, where we have used again equations \eqref{eq1} and \eqref{eq2}.
Thus, if the scalar curvature is constant, then $\Delta^\perp$ is totally geodesic and the desired result follows from 
Propositions \ref{prop} and \ref{conflat}. On the other hand, if the scalar is not constant then $\Delta^\perp$ is 
totally umbilical and the desired result follows from Propositions \ref{prop} and \ref{conflat2}. Finally, 
if $M^n$ is complete then the result is immediate and the proof is complete. \qed

\vspace{2ex}


\begin{thebibliography}{lll}





\bib{Blair}{article}{
   author={Blair, D. E.},
   title={On a generalization of the catenoid},
   journal={Canad. J. Math.},
   volume={27},
   date={1975},
   pages={231--236},
}







\bib{DCD}{article}{
   author={do Carmo, M.},
   author={Dajczer, M.},
   title={Rotation hypersurfaces in spaces of constant curvature},
   journal={Trans. Amer. Math. Soc.},
   volume={277},
   date={1983},
   number={2},
   pages={685--709},
}
		




\bib{cahe06}{article}{
   author={Carron, G.},
   author={Herzlich, M.},
   title={Conformally flat manifolds with nonnegative Ricci curvature},
   journal={Compos. Math.},
   volume={142},
   date={2006},
   number={3},
   pages={798--810},
}

\bib{car17}{article}{
   author={Cartan, E.},
   title={La d\'{e}formation des hypersurfaces dans l'espace conforme r\'{e}el \`a $n
   \ge 5$ dimensions},
   language={French},
   journal={Bull. Soc. Math. France},
   volume={45},
   date={1917},
   pages={57--121},
}




\bib{cadjnd11}{article}{
   author={Catino, G.},
   author={Djadli, Z.},
   author={Ndiaye, C. B.},
   title={A sphere theorem on locally conformally flat even-dimensional
   manifolds},
   journal={Manuscripta Math.},
   volume={136},
   date={2011},
   number={1-2},
   pages={237--247},
}

\bib{cat16}{article}{
   author={Catino, G.},
   title={On conformally flat manifolds with constant positive scalar
   curvature},
   journal={Proc. Amer. Math. Soc.},
   volume={144},
   date={2016},
   number={6},
   pages={2627--2634},
}




\bibitem{dajczer}  M. Dajczer,
{\it ``Submanifolds and Isometric Immersions"},
Math. Lecture Ser. 13, Publish or Perish Inc. Houston, 1990.



\bib{dt}{book}{
   author={Dajczer, M.},
   author={Tojeiro, R.},
   title={Submanifold theory beyond an introduction},
   series={Universitext},
   publisher={Springer US},
   date={2019},
}




\bib{df96}{article}{
   author={Dajczer, M.},
   author={Florit, L.A.},
   title={On conformally flat submanifolds},
   journal={Comm. Anal. Geom.},
   volume={4},
   date={1996},
   number={1-2},
   pages={261--284},
}

\bib{df99}{article}{
   author={Dajczer, M.},
   author={Florit, L.},
   title={Euclidean conformally flat submanifolds in codimension two
   obtained as intersections},
   journal={Proc. Amer. Math. Soc.},
   volume={127},
   date={1999},
   number={1},
   pages={265--269},
}


\bibitem{dov18}  M. Dajczer, C.-R. Onti and Th. Vlachos,
{\it ``Conformally flat submanifolds with flat normal bundle"}.
ArXiv e-prints (2018), available at {\sf https://arxiv.org/abs/1810.06968}.


\bib{mmf85}{article}{
   author={do Carmo, M.},
   author={Dajczer, M.},
   author={Mercuri, F.},
   title={Compact conformally flat hypersurfaces},
   journal={Trans. Amer. Math. Soc.},
   volume={288},
   date={1985},
   number={1},
   pages={189--203},
}

\bib{discala}{article}{
   author={Di Scala, A. J.},
   title={Minimal immersions of K\"{a}hler manifolds into Euclidean spaces},
   journal={Bull. London Math. Soc.},
   volume={35},
   date={2003},
   number={6},
   pages={825--827},
}


\bib{dote11}{article}{
   author={Donaldson, N.},
   author={Terng, C.-L.},
   title={Conformally flat submanifolds in spheres and integrable systems},
   journal={Tohoku Math. J. (2)},
   volume={63},
   date={2011},
   number={2},
   pages={277--302},
}

\bib{har70}{article}{
   author={Hartman, P.},
   title={On the isometric immersions in Euclidean space of manifolds with
   nonnegative sectional curvatures. II},
   journal={Trans. Amer. Math. Soc.},
   volume={147},
   date={1970},
   pages={529--540},
}


\bib{ku49}{article}{
   author={Kuiper, N. H.},
   title={On conformally-flat spaces in the large},
   journal={Ann. of Math. (2)},
   volume={50},
   date={1949},
   pages={916--924},
}

\bib{ku50}{article}{
   author={Kuiper, N. H.},
   title={On compact conformally Euclidean spaces of dimension $>2$},
   journal={Ann. of Math. (2)},
   volume={52},
   date={1950},
   pages={478--490},
}

\bib{mo77}{article}{
   author={Moore, J. D.},
   title={Conformally flat submanifolds of Euclidean space},
   journal={Math. Ann.},
   volume={225},
   date={1977},
   number={1},
   pages={89--97},
}



\bib{mm78}{article}{
   author={Moore, J. D.},
   author={Morvan, J.-M.},
   title={Sous-vari\'{e}t\'{e}s conform\'{e}ment plates de codimension quatre},
   language={French, with English summary},
   journal={C. R. Acad. Sci. Paris S\'{e}r. A-B},
   volume={287},
   date={1978},
   number={8},
   pages={A655--A657},
}

\bib{no90}{article}{
   author={N\"olker, S.},
   title={Isometric immersions with homothetical Gauss map},
   journal={Geom. Dedicata},
   volume={34},
   date={1990},
   number={3},
   pages={271--280},
}

\bib{no93}{article}{
   author={Noronha, M. H.},
   title={Some compact conformally flat manifolds with nonnegative scalar
   curvature},
   journal={Geom. Dedicata},
   volume={47},
   date={1993},
   number={3},
   pages={255--268},
}
		
\bib{pin85}{article}{
   author={Pinkall, U.},
   title={Compact conformally flat hypersurfaces},
   conference={
      title={Conformal geometry},
      address={Bonn},
      date={1985/1986},
   },
   book={
      series={Aspects Math., E12},
      publisher={Vieweg, Braunschweig},
   },
   date={1988},
   pages={217--236},
}

\bib{sy88}{article}{
   author={Schoen, R.},
   author={Yau, S.-T.},
   title={Conformally flat manifolds, Kleinian groups and scalar curvature},
   journal={Invent. Math.},
   volume={92},
   date={1988},
   number={1},
   pages={47--71},
}




\end{thebibliography}
\end{document}